\documentclass[a4paper,12pt]{article}
\usepackage{amssymb}
\usepackage{amsmath}
\usepackage{latexsym}
\usepackage{amsthm}
\usepackage{graphicx}
\usepackage{tikz}
\usetikzlibrary{svg.path}
\usepackage[colorlinks=true,linkcolor=black,citecolor=black,urlcolor=black]{hyperref}
\usepackage[top=1.5cm,left=1.5cm,bottom=1.5cm,right=1.5cm]{geometry}

\newtheorem{theorem}{Theorem}[section]

\newtheorem*{lemma*}{Lemma}

\newenvironment{config}{\vspace*{0.5ex}\tt\small\setlength{\parskip}{1pt}}{}

\date{}
\DeclareMathOperator{\STS}{STS}
\DeclareMathOperator{\KTS}{KTS}
\DeclareMathOperator{\PG}{PG}

\begin{document}
\title{Properties of Steiner triple systems of order 21}
\author{Grahame Erskine and Terry S. Griggs\\
        School of Mathematics and Statistics\\
        The Open University\\
        Walton Hall\\
        Milton Keynes MK7 6AA\\
        UNITED KINGDOM\\
        \texttt{grahame.erskine@open.ac.uk}\\
        \texttt{terry.griggs@open.ac.uk}}
\maketitle

\begin{abstract}
\noindent Properties of the 62,336,617 Steiner triple systems of order 21 with a non-trivial automorphism group are examined. In particular, there are 28 which have no parallel class, six that are 4-chromatic, five that are 3-balanced, 20 that avoid the mitre, 21 that avoid the crown, one that avoids the hexagon and two that avoid the prism. All systems contain the grid. None have a block intersection graph that is 3-existentially closed.
\end{abstract}


\section{Introduction}\label{sec:intro}
  A \emph{Steiner triple system} of \emph{order} $v$, denoted by $\STS(v)$, is an ordered pair $S=(V,\mathcal{B})$ where $V$ is a \emph{base set} of \emph{elements} or \emph{points} of cardinality $v$, and $\mathcal{B}$ is a collection of \emph{triples}, also called \emph{blocks}, which collectively have the property that every pair of distinct elements of $V$ is contained in precisely one triple. Such systems exist if and only if $v \equiv 1$ or $3 \pmod 6$; first proved by Kirkman in 1847~\cite{KIRK}. Two Steiner triple systems $(V,\mathcal{B})$ and $(W,\mathcal{D})$ are said to be \emph{isomorphic} if there exists a one-one mapping $\phi:V \to W$ such that every triple $B \in \mathcal{B}$ maps to a triple $\phi(B) \in \mathcal{D}$. To within isomorphism the $\STS(7)$ and $\STS(9)$ are unique. Already in the nineteenth century it was known that there are two non-isomorphic $\STS(13)$s~\cite{D},~\cite{Z}. The number of non-isomorphic $\STS(15)$s is 80; first determined by White, Cole and Cummings in 1919~\cite{WCC} and verified by Hall and Swift in 1955~\cite{HS} in one of the first uses of digital computers in Combinatorial Design Theory. Nearly fifty years then elapsed before Kaski and \"{O}sterg\r{a}rd~\cite{KO} enumerated the number of non-isomorphic $\STS(19)$s; there are 11,084,874,829 which are now stored in compact notation on the Internet~\cite{STS19}. They have been extensively analysed with the results appearing in~\cite{EIGHT}.

  Recently, Heinlein and \"{O}sterg\r{a}rd~\cite{HO} considered the next order for which Steiner triple systems exist and determined the number of non-isomorphic $\STS(21)$s: 14,796,207,517,873,771. It would be good to make a thorough analysis of these too, but the fact that there are more than a million times more systems than for $\STS(19)$ seems to make both storage and computation infeasible in the near future. However, previous to this enumeration, Kaski~\cite{K1} had determined all non-isomorphic $\STS(21)$s having a non-trivial automorphism group, by comparison a mere 62,336,617. It is feasible to investigate some of their properties; indeed Kaski himself already did so in the paper (automorphism groups, anti-Pasch systems). Two further papers~\cite{KO1},~\cite{KO2} also deal respectively with Kirkman triple systems and sparse systems. It is our aim to add to this analysis and identify systems which may be of particular interest. In the interests of economy of space, where such systems already appear in the literature we will not repeat them but instead point to their reference. Systems which have not already appeared in the literature will be listed. The focus will be mainly but not exclusively on parallel classes, colourings and configurations. For precise definitions of these concepts see later in the relevant sections. Our investigations also lead to a number of interesting questions. 
  
  When providing listings of particular systems, we will often employ the following compact notation which is common in the literature, see for example~\cite{FGG1}. The set $V$ of points is taken to be $\{0,1,2,\ldots,20\}$. The points $10,\ldots,20$ are represented by the letters $a,\ldots,k$ respectively. The 70 blocks of an $\STS(21)$ are represented by a string of symbols $s_1 s_2 \ldots s_{70}$. Using the usual lexicographical order, the symbol $s_i$ is the largest element $z_i$ in the $i$th triple $\{x_i,y_i,z_i\}$, where $x_i<y_i<z_i$. The remaining two elements implicitly have the property that there is no pair $x'_i<y'_i$ such that $\{x'_i,y'_i\}$ does not appear in an earlier triple, and either (i) $x'_i<x_i$ or (ii) $x'_i=x_i$ and $y'_i<y_i$.
  
\section{Parallel classes}\label{sec:parclass}
In an $\STS(v)$, $S=(V,\mathcal{B})$ where $v \equiv 3 \pmod 6$, a \emph{parallel class} or \emph{resolution class} is a set of blocks which contain every element precisely once. If the blocks of $\mathcal{B}$ can be partitioned into parallel classes then the $\STS(v)$ is said to be \emph{resolvable} and the set of parallel classes is called a \emph{resolution}. Such an $\STS(v)$ together with the resolution is called a \emph{Kirkman triple system} and denoted by $\KTS(v)$. The  $\STS(v)$ is said to be the \emph{underlying} Steiner triple system of the $\KTS(v)$. A Steiner triple system may underlie more than one Kirkman triple system. If an $\STS(v)$ underlies two $\KTS(v)$s and the latter have the property that every pair of parallel classes, one from each of the two systems, have either zero or one block in common then the $\STS(v)$ is said to be \emph{doubly resolvable}. Kirkman triple systems of order 21 with a non-trivial automorphism group were considered by Cohen, Colbourn, Ives and Ling~\cite{CCIL} and later by Kokkala and \"{O}sterg\r{a}rd~\cite{KO1}.  The latter paper further considered those systems whose underlying Steiner triple system has a non-trivial automorphism group. They obtained the following result.

\begin{theorem}\label{thm:underlie}\textnormal{(Kokkala and \"{O}sterg\r{a}rd)}\\
The Steiner triple systems of order 21 with a non-trivial automorphism group underlie 66,937 non-isomorphic Kirkman triple systems  with a non-trivial automorphism group and a further 1,992 Kirkman triple systems having only the identity automorphism. None of these systems is doubly resolvable.  
\end{theorem}

The reader is referred to the papers and the Tables therein for further information. In particular Table 8 of~\cite{CCIL} and Table 5 of~\cite{KO1} detail how the $\KTS(21)$ with an automorphism of order 3 and 2 respectively are distributed amongst the underlying $\STS(21)$. Table~\ref{tab:kts} shows the same information for ALL 68,929 $\KTS(21)$ having an underlying $\STS(21)$ with a non-trivial automorphism group and may be of interest. There are 55,900 such systems and Table~\ref{tab:numres} shows the number of resolutions that all of them have. One system stands out: that with 12,480 resolutions. Not surprisingly it is the one which underlies 18 Kirkman triple systems and also has 406 parallel classes, the most of any of the $\STS(21)$s. It is the direct product of the unique $\STS(7)$ and $\STS(3)$ systems and has automorphism group of order $168 \times 6 = 1008$. We represent the points of the system by $V = \mathbb{Z}_7 \times \{A, B, C \}$ and the blocks by $\mathcal{B} =\{ \{(i,X),(i+1,Y),(i+3,Z)\}, i \in \mathbb{Z}_7, \{X,Y,Z\}=\{A,B,C\}$ \textnormal{or} $X=Y=Z \in \{A,B,C\}, \{(i,A),(i,B),(i,C)\}, i \in \mathbb{Z}_7\}$. This system is isomorphic to the cyclic system on the points of $\mathbb{Z}_{21}$ generated by the blocks $\{0,1,5\}$, $\{0,2,10\}$, $\{0,3,9\}$ and $\{0,7,14\}$ under the mapping $i\mapsto i+1$: it is system C3 as listed in~\cite{MPR}.

\newpage
We remark that the system containing the next largest number of parallel classes, which is 294, is not resolvable. This system has an automorphism group of order 72 and is listed below.

\begin{config}
2468bcfgjk578cbgfkj867degfik7ihdekj8gfjkihaihkjghiedjkjkceihfgkjhikjhi
\end{config}

\begin{table}[h]
\begin{center}
 \begin{tabular}{|cc|cc|cc|}
 \hline 
\#KTS & \#STS & \#KTS & \#STS & \#KTS & \#STS\\
\hline
 1 & 45,604 & 7 & 9 & 13 & 1\\
 2 & 8,629 & 8 & 19 & 14 & 3\\
 3 & 1,025 & 9 & 5 & 15 & 1\\
 4 & 459 & 10 & 3 & 16 & 1\\
 5 & 99 & 11 & 3 & 18 & 1\\
 6 & 37 & 12 & 1 & ~ & ~\\
\hline 
\end{tabular}
\caption{Underlying $\STS(21)$s}
\label{tab:kts}
\end{center}
\end{table}

\begin{table}[h]
\begin{center}
 \begin{tabular}{|cc|cc|cc|}
 \hline 
\#Res & \#STS & \#Res & \#STS & \#Res & \#STS\\
\hline
 1 & 44,651 & 13 & 3 & 28 & 2\\
 2 & 9,208 & 14 & 3 & 32 & 4\\
 3 & 1,045 & 15 & 3 & 36 & 1\\
 4 & 639 & 16 & 5 & 44 & 2\\
 5 & 97 & 17 & 1 & 62 & 1\\
 6 & 90 & 18 & 3 & 64 & 1\\
 7 & 15 & 19 & 2 & 264 & 1\\
 8 & 75 & 20 & 3 & 336 & 1\\
 9 & 23 & 21 & 1 & 448 & 1\\
10 & 7 & 22 & 2 & 12,480 & 1\\
12 & 8 & 24 & 1 & ~ & ~\\
\hline
\end{tabular}
\caption{Numbers of resolutions}
\label{tab:numres}
\end{center}
\end{table}

Turning now to systems with no parallel classes, it is known that there are 12 non-isomorphic 4-rotational $\STS(21)$s, i.e. having an automorphism consisting of a fixed point and four 5-cycles, and there are given in~\cite{MR}. We have determined that there are a total of 28 non-isomorphic $\STS(21)$s with a non-trivial automorphism which have no parallel class. The other 16 seem not to be known. All have an automorphism of order 3; in 9 of these this acts as seven 3-cycles and the other 7 as three fixed points and six 3-cycles, We list the 16 new systems below.

The 9 systems with seven 3-cycles:

\begin{config}
fh6idgkacjge7jihdbkkc8fjbei5acfgjidbhgkj9efhikkjcehicdhedgikjkkjfhjgki

kj69cfedgiid7acghejbe8hdcfk5ekicgjicjdkhdjekhibagfk9hgkfhjkjfighkjiijk

269cikhdjgd7aijefkhbe8kjgciha9chkjibfdikjegjkia9jhgbkfhifgfhdjgekekjki

kj697acfhii87abdfgjb689ehgkgfehkjihcfikjdjgkiikcjfhjdkghiehggfkhijejkk

kj6f7dachii87gedbfjh68cbegkji9gdkfkaehgibfchjdcjhikekfjigijkhkfjgkjikh

26fdejcigke7gkchdjihc8dikfejbajgidk9hkjeiifkejcejgfdgkhfhiedcjhkhkjgik

kj69dcbegiie7a9dcjhbc8eadhkgf9ikjhhjaikfbkjgi8dfikegijchkjfhgkghjfkikj

kj6f7acdgii87gbedhjh689cefkjigfbkekhbgeifahdj8hkgjifijhkjgikfgkhjfkikj

kj6b7cafhii879dbfgja68ebhgkcehfgjkdfhgkjgfhikikjhcfjdkfggeihkijikjejik
\end{config}

\newpage
The 7 systems with three fixed points and six 3-cycles:

\begin{config}
2678cdeijk9abdecjkiijk9abfghgfbaekjhb9dika9jeikjifhigjhhkghgfjkfkgjhik

2678cdeijk9abfghijkdecjkighfcebjgkidk9hjiaifkjdcighehjffkgbkjkijiehhgk

2678cdeijk9abcdejkicdefghijkhgafjikfgbkij9hikj8kdijeijkjckibfhghhggfhk

2678cdeijk9abfghijk786kijhfgjibgekhk9hdigafjehedbgjkcbhkjafikgjikhjfik

2678cdeijk9abdecjkiijkab9hfgcejgfhkdhkfgiifhgjgfeckhdckedjkjiijkikgjhh

2678cdeijk9abdecjkiijk9abfgh5cajhgkbdkhfie9igfjhgfkjfgikhjkgfjihjkikeh

2678cdeijk9abdecjki867ijkhfgjidhgekkehfdicgfejjibkhgkaifhjbghkfgijhkjk
\end{config}

In addition there are 64 systems containing just a single parallel class.

\section{Colourings}\label{sec:colour}
A \emph{(weak) colouring} of an $\STS(v)$, $S=(V,\mathcal{B})$ is  a mapping $\phi:V \to C$ where $C$ is a set of cardinality $m$  whose elements are called \emph{colours}, such that $ |\phi(B)|  > 1$ for all $B \in \mathcal{B}$, i.e. no block is monochromatic.The system is said to be \emph{$m$-colourable}. The \emph{chromatic number} of $S$, $\chi(S)$ is the smallest value of $m$ for which $S$ admits a colouring with $m$ colours. We say that the Steiner triple system is \emph{$m$-chromatic}. All $\STS(v)$, $7 \leq v \leq 19$ are 3-chromatic. However for $v = 21$ , there exist 4-chromatic systems. The first of these was given by Haddad \cite{H} and five further ones were constructed in~\cite{FGG1} by applying various trades to Haddad's system. Respectively in the order in which they appear in~\cite{FGG1} their automorphism groups have order 108, 12, 12, 4, 3 and 18. Thus the system found by Haddad is the unique $\STS(21)$ with an automorphism group of order 108. We have determined that these six systems are the only 4-chromatic $\STS(21)$s having a non-trivial automorphism group. We have also tried various methods to try to construct further 4-chromatic $\STS(21)$s with no success. But with the large number of $\STS(21)$s it would be highly speculative and indeed rather foolish to conjecture that these are the only six 4-chromatic $\STS(21)$s.

With regard to the remaining 3-chromatic $\STS(21)$s, we focus on the cardinalities of the colour classes. An  $\STS(v)$ is said to be \emph{equitably $m$-chromatic} if it is $m$-chromatic and the cardinalities of the colour classes differ by at most one. It is \emph{$m$-balanced} if every $m$-colouring is equitable and \emph{$m$-unbalanced} if every $m$-colouring is not equitable. It was proved in~\cite{FGG1} that every 3-chromatic $\STS(21)$ has a 3-colouring in which the cardinalities of the colour classes are either $(8,7,6)$ or $(7,7,7)$.
Our investigations show that in fact every 3-chromatic $\STS(21)$ with a non-trivial automorphism group is equitably 3-chromatic. Further all but five also have a 3-colouring with colour classes with cardinalities $(8,7,6)$. Two of these systems appear in~\cite{FGG1} but the other three were previously unknown. For completeness we list all five of these 3-balanced $\STS(21)$s below.

\begin{config}
2468bcfgjk56jkcbgfi65ihacegkdahjgik9eifkhjgbekijhcfjdhkifekgdjkijhikhj

2678cdeijk9abfghjkiijkecdfghkjdchgfidehfgcegfh8ikjhkjihjikgedcekijikjk

2468bcfgjk578dehijk867jkchig7ihfgkj8kjihgfcbihkjadeihkgfjkeiedgfkjhikj

2468bcfgjk56jkcbgfi65ihacegkdahkgij9eifjhkgbekijhcfjdhkifekgdjkijhijhk

2468bcfgjk578achijk867cbjkhi7ihdekj8fgkjihhijkgfjkhiegedgfikgfedkjkjih
\end{config}

Let  $v \equiv 3 \pmod 6$. A simple counting argument shows that in any equitable colouring of a 3-chromatic $\STS(v)$, the number of blocks which contain a point from each of the three colour classes is $v/3$. Call this set of blocks the \emph{rainbow set}. This is suggestive that the rainbow set might be a parallel class. We checked this for the 3-chromatic $\STS(21)$s having a non-trivial automorphism group and found that all have an equitable 3-colouring in which the rainbow set is NOT a parallel class. Only 465,006, less than 1\% of the total, have such a colouring where the rainbow set is a parallel class. No system had an equitable 3-colouring based solely on the rainbow set being a parallel class and it is an interesting question of whether such a system exists among those with only the identity automorphism. Certainly such systems exist for other orders. The $\STS(15)$ formed by the point-line design of the projective geometry $\PG(3,2)$ has up to isomorphism a unique 3-colouring which is equitable and the rainbow set is a parallel class~\cite{P,FHW}. Further in~\cite{F}, Forbes lists a number of systems which are uniquely equitably 3-colourable.In those of order 39 and 57, the rainbow set is a parallel class. 

\section{Configurations}\label{sec:config}
A $(k,\ell)$-\emph{configuration} in a Steiner triple system is a collection of $\ell$ blocks which contain $k$ points. It is said to be \emph{even} if every point occurs in an even number of blocks. It is \emph{full} if no point occurs in just one block. Configurations with $k=\ell+2$ play a special role in the structural theory of $\STS(v)$. Apart from the single block, the smallest example is the so-called \emph{Pasch configuration}, the unique configuration with $(k,\ell)=(6,4)$. The four blocks are isomorphic to $\{x,y,z\}$, $\{x,b,c\}$, $\{a,y,c\}$, $\{a,b,z\}$. There are two $(7,5)$-configurations. One of them, the \emph{mia} is obtained by extending the Pasch configuration with one extra block through any pair of uncovered points and so is not full. The other is the \emph{mitre}; five blocks isomorphic to $\{x,a,d\}$, $\{x,b,e\}$, $\{x,c,f\}$, $\{a,b,cz\}$, $\{d,e,f\}$.
There are five $(8,6)$-configurations (and one $(7,6)$-configuration) but three are obtained  by extending the mia or the mitre and are not full. The other two (full) configurations are the \emph{6-cycle} or \emph{hexagon}, six blocks isomorphic to 
$\{x,a,b\}$, $\{x,c,d\}$, $\{x,e,f\}$, $\{y,a,f\}$, $\{y,b,c\}$, $\{y,d,e\}$.
and the \emph{crown}, six blocks isomorphic to $\{x,y,z\}$, $\{x,a,b\}$, $\{y,c,d\}$, $\{z,a,c\}$, $\{w,x,d\}$, $\{w,y,b\}$.  
Steiner triple systems which do not contain some or all of these configurations are of particular interest. For $n \geq 4$, an $\STS(v)$ with the property that it contains no $(\ell+2,\ell)$-configuration for $4 \leq \ell \leq n$ is said to be $n$-\emph{sparse}.

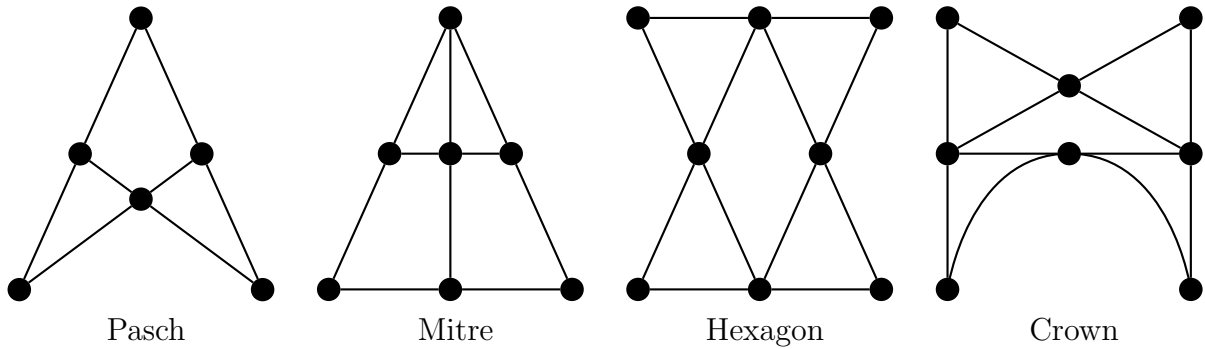
\begin{figure}\centering
    \begin{tabular}{cccc}
        \begin{tikzpicture}[x=0.2mm,y=0.15mm,thick,vertex/.style={circle,draw,minimum size=8,inner sep=0,fill=black}]
	\node at (80,-93.3) [vertex] (v1) {};
	\node at (-80,-93.3) [vertex] (v2) {};
	\node at (0,-13.3) [vertex] (v3) {};
	\node at (40,26.7) [vertex] (v4) {};
	\node at (-40,26.7) [vertex] (v5) {};
	\node at (0,146.7) [vertex] (v6) {};
	\draw (v5) to (v6);
	\draw (v2) to (v5);
	\draw (v4) to (v6);
	\draw (v1) to (v4);
	\draw (v2) to (v3);
	\draw (v3) to (v4);
	\draw (v1) to (v3);
	\draw (v3) to (v5);
\end{tikzpicture} &
        \begin{tikzpicture}[x=0.2mm,y=0.15mm,thick,vertex/.style={circle,draw,minimum size=8,inner sep=0,fill=black}]
	\node at (90,-90) [vertex] (v1) {};
	\node at (-70,-90) [vertex] (v2) {};
	\node at (10,30) [vertex] (v3) {};
	\node at (50,30) [vertex] (v4) {};
	\node at (-30,30) [vertex] (v5) {};
	\node at (10,150) [vertex] (v6) {};
	\node at (10,-90) [vertex] (v7) {};
	\draw (v5) to (v6);
	\draw (v2) to (v5);
	\draw (v4) to (v6);
	\draw (v1) to (v4);
	\draw (v3) to (v4);
	\draw (v3) to (v5);
	\draw (v2) to (v7);
	\draw (v1) to (v7);
	\draw (v3) to (v6);
	\draw (v3) to (v7);
\end{tikzpicture} &
        \begin{tikzpicture}[x=0.2mm,y=0.15mm,thick,vertex/.style={circle,draw,minimum size=8,inner sep=0,fill=black}]
	\node at (90,-90) [vertex] (v1) {};
	\node at (-70,-90) [vertex] (v2) {};
	\node at (50,30) [vertex] (v4) {};
	\node at (-30,30) [vertex] (v5) {};
	\node at (10,150) [vertex] (v6) {};
	\node at (10,-90) [vertex] (v7) {};
	\node at (-70,150) [vertex] (v8) {};
	\node at (90,150) [vertex] (v9) {};
	\draw (v2) to (v7);
	\draw (v1) to (v7);
	\draw (v6) to (v8);
	\draw (v6) to (v9);
	\draw (v5) to (v8);
	\draw (v5) to (v7);
	\draw (v2) to (v5);
	\draw (v5) to (v6);
	\draw (v4) to (v6);
	\draw (v1) to (v4);
	\draw (v4) to (v9);
	\draw (v4) to (v7);
\end{tikzpicture} &
        \begin{tikzpicture}[x=0.2mm,y=0.15mm,thick,vertex/.style={circle,draw,minimum size=8,inner sep=0,fill=black}]
	\node at (110,-90) [vertex] (v1) {};
	\node at (-50,-90) [vertex] (v2) {};
	\node at (30,30) [vertex] (v3) {};
	\node at (110,30) [vertex] (v4) {};
	\node at (-50,30) [vertex] (v5) {};
	\node at (-50,150) [vertex] (v6) {};
	\node at (30,90) [vertex] (v7) {};
	\node at (110,150) [vertex] (v8) {};
	\draw (v2) to (v5);
	\draw (v1) to (v4);
	\draw (v3) to (v4);
	\draw (v3) to (v5);
	\draw (v6) to (v7);
	\draw (v7) to (v8);
	\draw (v4) to (v7);
	\draw (v5) to (v7);
	\draw (v4) to (v8);
	\draw (v5) to (v6);
	\draw plot [smooth,tension=1.5] coordinates {(v2) (v3) (v1)};
\end{tikzpicture} \\
        Pasch & Mitre & Hexagon & Crown
    \end{tabular}
    \caption{The full $(\ell+2,\ell)$ configurations on at most six blocks}
    \label{fig:config}
\end{figure}

\begin{figure}\centering
	\begin{tabular}{cc}
		\begin{tikzpicture}[x=0.2mm,y=0.2mm,thick,vertex/.style={circle,draw,minimum size=8,inner sep=0,fill=black}]
	\node at (110,-10) [vertex] (v1) {};
	\node at (-50,-10) [vertex] (v2) {};
	\node at (30,70) [vertex] (v3) {};
	\node at (110,70) [vertex] (v4) {};
	\node at (-50,70) [vertex] (v5) {};
	\node at (-50,150) [vertex] (v6) {};
	\node at (30,150) [vertex] (v7) {};
	\node at (110,150) [vertex] (v8) {};
	\node at (30,-10) [vertex] (v9) {};
	\draw (v2) to (v5);
	\draw (v1) to (v4);
	\draw (v3) to (v4);
	\draw (v3) to (v5);
	\draw (v6) to (v7);
	\draw (v7) to (v8);
	\draw (v4) to (v8);
	\draw (v5) to (v6);
	\draw (v2) to (v9);
	\draw (v1) to (v9);
	\draw (v3) to (v9);
	\draw (v3) to (v7);
\end{tikzpicture} &
		\begin{tikzpicture}[x=0.16mm,y=0.16mm,thick,vertex/.style={circle,draw,minimum size=8,inner sep=0,fill=black}]
	\node at (46.9,98.9) [vertex] (v1) {};
	\node at (100.5,82.8) [vertex] (v2) {};
	\node at (190,10) [vertex] (v3) {};
	\node at (-50,10) [vertex] (v4) {};
	\node at (70,210) [vertex] (v5) {};
	\node at (130,110) [vertex] (v6) {};
	\node at (70,10) [vertex] (v7) {};
	\node at (10,110) [vertex] (v8) {};
	\node at (61.6,43.3) [vertex] (v9) {};
	\draw (v3) to (v6);
	\draw (v5) to (v6);
	\draw (v3) to (v7);
	\draw (v4) to (v7);
	\draw (v5) to (v8);
	\draw (v4) to (v8);
	\draw (v1) to (v9);
	\draw (v1) to (v2);
	\draw (v2) to (v9);
	\draw (v7) to (v9);
	\draw (v2) to (v6);
	\draw (v1) to (v8);
\end{tikzpicture} \\
		Grid & Prism 
	\end{tabular}
	\caption{The even $(9,6)$-configurations}
	\label{fig:evenconfig}
\end{figure}
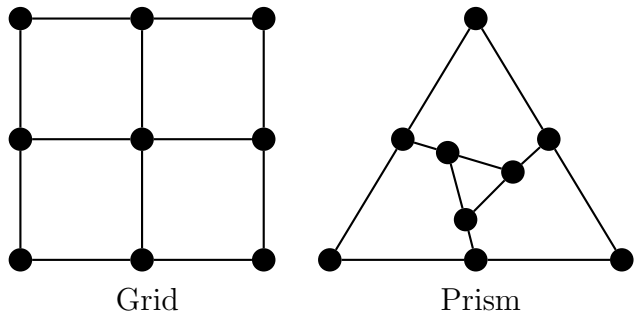

So a 4-sparse $\STS(v)$ is just one which contains no Pasch configurations. Such systems are more commonly known as \emph{anti-Pasch} and they exist for all $v \equiv 1$ or $3 \pmod 6$, $v \neq 7, 13$ \cite {LCGG}, \cite {GGW}. Analogously a system which contains no mitres is called \emph{anti-mitre} and a 5-sparse $\STS(v)$ is one which avoids both the Pasch and the mitre configurations. Anti-mitre $\STS(v)$ exist for all $v \equiv 1$ or $3 \pmod 6$, $v \neq 9$ but the exact spectrum for 5-sparse $\STS(v)$ is not completely determined, \cite{CMRS}, \cite{L}, \cite{F1}, \cite{F2}, \cite{W1}, \cite{W2}.
Kokkala and  \"{O}sterg\r{a}rd \cite{KO2} have enumerated all anti-Pasch $\STS(21)$s; there are 83,003,869 non-isomorphic such systems, of which just three are 5-sparse and none are 6-sparse. See \cite{KOZ} for a listing.
Previous to this, Kaski \cite{K1} had determined that there are 958 non-isomorphic anti-Pasch $\STS(21)$s with a non-trivial automorphism group.
Just one of these is cyclic: it is system C2 as listed in~\cite{MPR}. On the base set $V = \mathbb{Z}_{21}$ it can be obtained by the action of the mapping $i \mapsto i+1$ on the blocks $\{0,1,3\}$, $\{0,4,12\}$, $\{0,5,15\}$, $\{0,7,14\}$.

On the same data set, our investigations show that there are 20 non-isomorphic anti-mitre systems including the three that are 5-sparse and 21 systems which do not contain the crown including one of those which is anti-mitre. This system which avoids both the mitre and the crown is cyclic and on the base set $V = \mathbb{Z}_{21}$ is obtained by the action of the mapping $i \mapsto i+1$ on the blocks $\{0,1,5\}$, $\{0,2,13\}$, $\{0,3,9\}$, $\{0,7,14\}$. It is system C5 as listed in~\cite{MPR}; its full automorphism group has order 882. We return to this system later in the paper.
There is also a unique system containing no hexagon, listed below:

\begin{config}
2468bcfgjkde78afgik78hijkcegjahegki9kidjfhgbdikfcehjgjdhkkfiejihhikjkj
\end{config}

Turning now to systems containing the largest numbers of these configurations, the $\STS(21)$ with a non-trivial automorphism group having the largest number of Pasch configurations is the unique system with automorphism group of order 108 listed below; it contains 117 Pasch configurations. This system has one fixed point (0) and three further point orbits of sizes 9, 9 and 2 respectively.

\begin{config}
2468bcfgjk578achijk867cbjkhi7ihfgkj8dekjihhijkgfjkhiegedgfkjgfedikkjih
\end{config}

The largest number of mitres in any such system is 252, occurring in 12 systems, Remarkably, all also contain precisely 252 hexagons and no crowns.
Two of these systems are cyclic: they are systems C1 and C3 as listed in~\cite{MPR}. On the base set $V = \mathbb{Z}_{21}$ they can be obtained by the action of the mapping $i \mapsto i+1$ on the blocks $\{0,1,3\}$, $\{0,4,12\}$, $\{0,5,11\}$, $\{0,7,14\}$ and $\{0,1,5\}$, $\{0,2,10\}$, $\{0,3,9\}$, $\{0,7,14\}$ respectively. Of the remainder, one is the unique Pasch-maximal system described above with an automorphism group of order 108. The others have automorphism groups varying in order from 6 to 54; none is point-transitive.

There is a unique $\STS(21)$ with a non-trivial automorphism group having the largest number of crowns. It contains 396 crowns, has automorphism  group of order 18 and is listed below. 

\begin{config}
2468bcfgjkbcdejkagi9afghicekh7fjigk8igkhjfaekgj9dfjkcdikbehjijhhikjikh
\end{config}

Possibly of greatest interest though is the unique system having the greatest number of hexagons, 441 in total. It is the same system which has no mitres or crowns listed above.

Turning now to even configurations, the smallest (least number of blocks) is obviously the Pasch configuration. The next smallest are two $(9,6)$-configurations, the \emph{grid} and the \emph{prism}. Respectively the blocks are isomorphic to $\{a,b,c\}$, $\{\ell,m,n\}$, $\{x,y,z\}$, $\{a,\ell,x\}$, $\{b,m,y\}$, $\{c,n,z\}$ and $\{a,b,n\}$, $\{a,c,m\}$, $\{b,c,\ell\}$, $\{\ell,y,z\}$, $\{m.x.z\}$, $\{n,x,y\}$. In \cite{FC}, Fujiwara and Colbourn 
proved that every Steiner triple system contains an even configuration of cardinality 8 or less. It follows therefore that systems avoiding either the grid of prism are of interest. However little seems to be known. The unique $\STS(7)$ avoids both configurations trivially because it does not contain enough points and the unique $\STS(9)$ avoids the prism but not the grid. It is readily verified by computer that both $\STS(13)$s and all 80 $\STS(15)$s contain both grids and prisms.  In \cite {EIGHT}, this was extended to show that the same is true for all $\STS(19)$s.

Investigating the $\STS(21)$s with a non-trivial automorphism group, we can report that there is no system which avoids the grid. However there are two systems which avoid the prism. One of these is the cyclic system C5 and the other is the unique system with automorphism group of order 294. It is also one of the 21 systems which do not contain the crown. It is listed below.

\begin{config}
2468acegik56ijbghfk659fkdjihcehdgkjkhfcijggajkehidbeifjkbhjcikkijhfjgk
\end{config}

\newpage
With regard to systems containing most grids and prisms, the largest number of grids in an $\STS(21)$ with a non-trivial automorphism group is 798 which occurs in two systems, both of which are cyclic (C1 and C3). The largest number of prisms is 1773, occurring in the unique system with automorphism group of order 54 and listed below.

\begin{config}
2468acegik857kfehij678gjidhk7bcfgjk8idkfjhdejkhichgkfjhifgjkebjgikjkhi
\end{config}

\section{Pasch trades}\label{sec:trades}
The Pasch configuration is a \emph{trade}, i.e. it can be replaced in a Steiner triple system by a different collection of blocks containing the same pairs, possibly to create a non-isomorphic Steiner triple system. Specifically the blocks $\{x,y,z\}$, $\{x,b,c\}$, $\{a,y,c\}$, $\{a,b,z\}$ are replaced by the blocks $\{a,b,c\}$, $\{a,y,z\}$, $\{x,c,z\}$, $\{x,y,c\}$.  It is known that any one of the 79 $\STS(15)$s which contain a Pasch configuration can be transformed to any other by a sequence of such Pasch trades \cite{G}.

This led the authors of~\cite{GGM1} to define the concept of \emph{twin Steiner triple systems} which are two $\STS(v)$s, each of which contains precisely one Pasch configuration which when switched produces the other system. Each twin of a pair of twin Steiner triple systems has the same automorphism group. However the twins themselves need not be isomorphic.
When they are, appropriately  they are called \emph{identical twins}. There are 838 non-isomorphic Steiner triple systems with a non-trivial automorphism group which contain one Pasch configuration, all in fact having an automorphism of order 3. This includes six pairs of twins which are listed below but no identical twins. There does however exist a pair of identical twin $\STS(21)$s, see~\cite {GGM1}, necessarily of course having only the identity automorphism. There are no systems containing precisely one hexagon (which also is a trade) or one crown, and just one system containing one mitre, given below.

The six pairs of twin systems:

\begin{config}
kj69cfegihid7acgjhfbe8hdfkggfaihjkhjbfik9kgjiedbikcbjiajkkeieijdjkhghk

269cfegihkd7acgjhfkbe8hdfkgjgfaihjkhjbfik9kgjiedbikcbjiajkkeieijdjkhgh

\vspace*{1ex}
kj69cfbdhiid7a9gefjbe8hacgka9ikjfhbijgkhkjhigbaefj9dgkheihgikjfijkkjik

269cfbdhikd7a9gefjkbe8hacgkja9ikjfhbijgkhkjhigbaefj9dgkheihgikjfijkkji

\vspace*{1ex}
ed69chbjikcd7a9fkjibe8gaijkji8cfkhk8dgih7hejgejfhkkcgfjdighkgfjhkiejki

269chbjeikd7a9fkejibe8gaidjkji8cfkhk8dgih7hejgejfhkkcgfjdighkgfjhkijki

\vspace*{1ex}
26bcgfakjid79ghbkijae8fhbjik87gehik6chfjkdfgkjkjiehdifjcedgkekjikjighh

a96bcgfkjibd79ghkijae8fhjik87gehik6chfjkdfgkjkjiehdifjcedgkebkjikjighh

\vspace*{1ex}
26bcgfakjid79ghbkijae8fhbjik87gcfjk6dhgkjefhikikjchejfkdeegidikjikjhgh

a96bcgfkjibd79ghkijae8fhjik87gcfjk6dhgkjefhikikjchejfkdeegidbikjikjhgh

\vspace*{1ex}
26h7jacfik87fkdbjgig68ibehkjikafdejjbegdk9cheidcjihgejkfhikghgkhifjkjk

a96h7jcfikb87fkdjgig68iehkjikafdejjbegdk9cheidcjihgejkfhikghbgkhifjkjk
\end{config}

The unique system containing a single mitre:

\begin{config}
2678cdeijk9abfghijkdeckijghfjibakfhkb9giha9hjg8jhgikjfhifkgbekiejidkjk
\end{config}

The authors of~\cite{EIGHT} extended this analysis. They considered Steiner triple systems containing precisely two Pasch configurations, say \textit{P} and \textit{Q}, such that when \textit{P} (respectively \textit{Q}) is switched what is obtained is a Steiner triple system  containing just one Pasch configuration \textit{P'} (respectively \textit{Q'}). Clearly \textit{P} and \textit{Q} must have a single common block. (They cannot have two common blocks; otherwise the Steiner triple system would contain at least one further Pasch configuration.) Considering only the $\STS(21)$s with a non-trivial automorphism group having precisely two Pasch configurations, there are four systems which have this property. They are listed below. In all four cases, the system admits an automorphism which exchanges its two Pasch configurations \textit{P} and \textit{Q}, and thus the two switched systems containing a single Pasch configuration are isomorphic. 

\begin{config}
2468acegik69beghjikcfa8jdehkdgkfbhjjaheigkicgkhcbgijkkdfijhfejkfijkhij

2468acegik69ejcfhgkgjh8kicdfchkeifj7bidgkh9jgkiadjikhgiefbhgjkfjhkjikj

2468acegikh8adjcfgk9igkcdfejd7egkijkbdhjgff9ijhijefheigkcjhkdkihkhkijj

2468acegikajhdf9gekgekihbcjfeb9hjki8gcdfikafdijjfkcheikdkgjhjikihjghjk
\end{config}

This situation does not occur with the $\STS(19)$s where there are nine such systems, but in every case the two systems containing a single Pasch configuration are non-isomorphic; see~\cite{EIGHT}. For the $\STS(21)$s, a system containing two Pasch configurations which when either is switched gives a system with one Pasch configuration and the two systems are not isomorphic is given below. The systems of course only have the identity automorphism.  

\begin{config}
2468acegikgkfbh9jdich8djegfkajkhdifcfeaigjbgijkighkecfgdjkikjeihjhkkhj    
\end{config}

\section{Other properties}\label{sec:other}
Let $(V,\mathcal{B})$ be an $\STS(v)$. For each pair $x,y \in V$, define a graph $G_{x,y}$ with vertex set $G \setminus \{x,y,z\}$
where $\{x,y,z\} \in \mathcal{B}$, and two vertices $u,v$ are joined by an edge if either $\{x,u,v\} \in \mathcal{B}$ or $\{y,u,v\} \in \mathcal{B}$. The graph $G_{x,y}$ is a union of cycles of even length greater than 2, and these can be recorded  as a list of cycle lengths in non-decreasing order. The \emph{cycle structure} of the Steiner triple system is the collection of all such \emph{cycle lists}. Within this framework, the greatest interest is in Steiner triple systems in which all cycle lists are the same; such systems are called \emph{uniform}. Of even greater interest are uniform Steiner triple systems in which each cycle list is $v-3$; such systems are called \emph{perfect} and only 14 of these are known \cite{GGM2}, \cite{FGG2}. Kaski \cite{K2}, \cite{KMO} has determined that there is no perfect $\STS(21)$. Amongst the $\STS(21)$s with a non-trivial automorphism group, we have determined that there is also no uniform system. Indeed all systems except one have at least three different cycle structures in their cycle lists.  Not surprisingly the exception is the system C5 as listed in~\cite{MPR} which contains the largest number of hexagons, There are 147 pairs with cycle structure $6, 6, 6$ and 63 pairs with cycle structure $4, 14$. 

An \emph{independent set} $I \subset V$ in a Steiner triple system $S=(V,\mathcal{B})$ is a set of points with the property that no block of $\mathcal{B}$ is contained in $I$. A \emph{maximum independent set} is an independent set of maximum cardinality. It is know that a maximum independent set of an $\STS(21)$ has cardinality 8, 9 or 10 \cite{FGG2}. The number of non-isomorphic $\STS(21)$s with a non-trivial automorphism group having each of these values is 15,614,086, 43,050,614 and 3,671,917 respectively.

The \emph{block intersection graph} of a Steiner triple system has one vertex for each block and an edge between two vertices when the corresponding blocks intersect. A graph $G=(W,\mathcal{E})$ is \emph{n-existentially closed} if for every $n$-element subset $S \subseteq W$ of vertices and for every subset $T \subseteq S$, there exists a vertex $x \notin S$ that is adjacent to every vertex in $T$ and not adjacent to every vertex in $S \setminus T$. The block intersection graph of an $\STS(v)$ is 2-existentially closed if and only if $v \geq 13$ \cite {FGG3}, \cite{MP}. It cannot be 3-existentially closed for $v \leq 15$ or $v \geq 25$. There are two $\STS(19)$s which possess 3-existentially closed block intersection graphs \cite{EIGHT}. So the only order for which this property is not determined is $v = 21$. We can report that no $\STS(21)$ with a non-trivial automorphism group has a block intersection graph which is 3-existentially closed.

\section{Extremal systems}\label{sec:extremal}
Much of this paper has been concerned with identifying $\STS(21)$s with a non-trivial automorphism group which are extremal with respect to a certain property. As we have shown, and perhaps not surprisingly, some systems are extremal with respect to a number of properties. We collect these systems together in this section which we hope will provide a succinct summary of the main results, presented by system rather than property. There are seven systems in all, four of which are cyclic.

Unique system with automorphism group of order 54.
This is the unique system containing the largest number of prisms (1773) and one of the 12 systems containing the largest number of mitres (252). 

Unique system with automorphism group of order 108.
This is the unique system containing the largest number of Pasch configurations (117) and one of the 12 systems containing the largest number of mitres (252). It has no crowns.

Unique system with automorphism group of order 294.
This is one of the two systems which contains no prisms. It also has no crowns. 

Cyclic system C1 with automorphism group of order 504.
This is one of the two systems which contains the largest number of grids (798) and one of the two \emph{cyclic} systems which contains the largest number of mitres (252). It has no crowns.  

Cyclic system C2 with automorphism group of order 21.
This is the unique \emph{cyclic} anti-Pasch $\STS(21)$.  

Cyclic system C3 with automorphism group of order 1008.
This is the other system which contains the largest number of grids (798), the other \emph{cyclic} system which contains the largest number of mitres (252) and no crowns. It is also the direct product of the Steiner triple systems $\STS(7)$ and $\STS(3)$ and has the largest number of resolutions (12,480) and parallel classes (406) and underlies the most Kirkman triple systems (18).
  
Cyclic system C5 with automorphism group of order 882.
This system is extremal with respect to a number of properties. It is the unique system which avoids both the mitre and the crown and is also the other system which contains no prisms. It is the unique system containing the largest number of hexagons (441) and the only system with just two cycle structures.

\section{Further questions}\label{sec:questions}
The analysis of the $\STS(21)$s with a non-trivial automorphism group presented in this paper naturally raises a number of questions relating to ALL $\STS(21)$s which would be answered if and when a complete investigation of these systems becomes possible. Some however might be resolved by alternative means, though we do not underestimate the tedium or difficulty of doing so. The first question comes from Section \ref{sec:parclass} and is simply stated.

\textbf{Question 1.} Does there exist a doubly resolvable $\STS(21)$?

A number of questions are raised from Section \ref{sec:colour}. The chromatic number of all $\STS(21)$s with a non-trivial automorphism group is 3, except for six systems which are 4-chromatic. It would be good to know whether these are the only systems with this property.

\textbf{Question 2.} Does there exist a 4-chromatic $\STS(21)$ with a trivial automorphism group?

In relation to 3-chromatic systems, we have noted that all $\STS(21)$s have a colouring in which the cardinalities of the colour classes are either $(8,7,6)$ or $(7,7,7)$. Amongst the $\STS(21)$s with a non-trivial automorphism group, all have a colouring with colour classes of cardinalities $(7,7,7)$, i.e. they are equitably 3-chromatic, and all but five also have a colouring with colour classes of cardinalities $(8,7,6)$. These five systems having only an equitable 3-colouring are 3-balanced. This leads to the next two questions.

\textbf{Question 3.} Do there exist further 3-balanced $\STS(21)$s?

\textbf{Question 4.} Is every 3-chromatic $\STS(21)$ also equitably 3-chromatic?

Also in Section \ref{sec:colour}, in an equitable colouring of a 3-chromatic $\STS(v)$, where $v \equiv 3 \pmod 6$, we introduced the concept of a rainbow set. This has the same cardinality ($v/3$) as that of a parallel class. However none of the $\STS(21)$s under consideration had a colouring in which the rainbow set is a parallel class and leads to the final question from the section.

\textbf{Question 5.} Does there exist a 3-chromatic $\STS(21)$ which has an equitable colouring in which the rainbow set is a parallel class? 

The main questions to arise from Section \ref{sec:config} concern the numbers of systems which avoid the mitre, crown, hexagon, grid or prism or any of these in combination. As we have noted, there are 83,003,689 anti-Pasch $\STS(21)$s of which just 958 have a non-trivial automorphism group. The number of anti-mitre, anti-crown and anti-hexagon systems with a non-trivial automorphism group are 20, 21 and 1 respectively which may indicate that there are not too many of such systems in total. But the most interesting problems concern the grid and the prism.

\textbf{Question 6.} Does there exist an anti-grid $\STS(21)$?

\textbf{Question 7.} Are there any anti-prism $\STS(21)$s other than the two identified in this paper?

Referring to Section  \ref{sec:trades}, we recall that Pasch trades partition the set of non-isomorphic $\STS(21)$s into equivalence classes. Based on admittedly very limited evidence from the situations for the $\STS(15)$s and $\STS(19)$s, the expectation is probably that most systems will be in one gigantic class (the \emph{peloton} to borrow a term from cycling) with some smaller classes such as twins and other small sets of systems. But this may not be the case; the peloton may be divided and the smaller classes may in fact contain thousands or even millions of systems, given the huge number of non-isomorphic $\STS(21)$s. A complete determination of these equivalence classes seems to be totally infeasible but any further information would undoubtedly be welcome.

Finally from section \ref{sec:other}, the main question concerns whether there exists an $\STS(21)$ whose block intersection graph is 3-existentially closed? As we noted $v=21$ is the only order for which this is not determined and we have shown that the answer is in the negative for $\STS(21)$s with a non-trivial automorphism group. So our final question is the following.

\textbf{Question 8,} Does there exist an $\STS(21)$, necessarily having only the identity automorphism, whose block intersection graph is 3-existentially closed?

Given the total number of $\STS(21)$s, we have been reluctant to make conjectures concerning the above questions but for this final one we would be surprised if such a system exists.

\section*{Acknowledgements}
We thank Petteri Kaski for making available to us the systems having non-trivial automorphism group identified in~\cite{K1}, and Patric \"{O}sterg\r{a}rd for helpful discussions during the drafting of this paper.


\end{document}